\documentclass{amsproc}%
\usepackage{graphicx}
\usepackage{amscd}
\usepackage{amsmath}%
\usepackage{amsfonts}%
\usepackage{amssymb}
\theoremstyle{plain}
\newtheorem{theorem}{Theorem}[section]

\newtheorem{corollary}[theorem]{Corollary}

\newtheorem{lemma}[theorem]{Lemma}

\theoremstyle{definition}
\newtheorem{definition}[theorem]{Definition}
\theoremstyle{remark}
\newtheorem*{acknowledgements}{Acknowledgements}
\numberwithin{equation}{section}

\begin{document}
\title[Cuntz algebras, loop groups and wavelets]{Representations of Cuntz algebras, loop groups and wavelets}
\author{Palle E. T. Jorgensen}
\address{Department of Mathematics\\
The University of Iowa\\
14 MacLean Hall\\
Iowa City, IA 52242-1419\\
U.S.A.}
\email{jorgen@math.uiowa.edu}
\thanks{Work supported in part by the U.S. National Science Foundation.}
\thanks{This paper represents an expanded version of an invited lecture given by the
author at the International Congress on Mathematical Physics, July 2000 in London.}

\begin{abstract}
A theorem of Glimm states that representation theory of an NGCR $C^{\ast}%
$-algebra is always intractable, and the Cuntz algebra $\mathcal{O}_{N}$ is a
case in point. The equivalence classes of irreducible representations under
unitary equivalence cannot be captured with a Borel cross section.
Nonetheless, we prove here that wavelet representations correspond to
equivalence classes of irreducible representations of $\mathcal{O}_{N}$, and
they are effectively labeled by elements of the loop group, i.e., the group of
measurable functions $A:\mathbb{T}\rightarrow\mathrm{U}_{N}(\mathbb{C})$.
These representations of $\mathcal{O}_{N}$ are constructed here from an orbit
picture analysis of the infinite-dimensional loop group.

\end{abstract}
\maketitle

\section{Introduction\label{section1}}

Recall the Cuntz algebra with $N$ generators $s_{0},\dots,s_{N-1}$ is the
$C^{\ast}$-algebra $\mathcal{O}_{N}$ on the relations%
\begin{equation}
s_{i}^{\ast}s_{j}=\delta_{ij}1\text{,\qquad}\sum_{i=0}^{N-1}s_{i}s_{i}^{\ast
}=1. \label{eq1.1}%
\end{equation}
Cuntz \cite{Cun77} showed that it is simple and infinite. (We will consider
$N$ finite only, $N=2,3,\dots$.) By a theorem of \cite{Gli60} there does not
exist a Borel section parameterizing the irreducible representations of
$\mathcal{O}_{N}$. Hence we shall restrict to special representations and
consider the questions of irreducibility and decomposition. In particular we
show that the elements in the loop group, i.e., all measurable maps%
\[
A:\mathbb{T}\longrightarrow\mathrm{U}_{N}(\mathbb{C}),
\]
parameterize the equivalence classes of wavelet representations. The wavelets
in turn correspond $1$-to-$1$ to representations
\[
\pi^{(A)}\in\operatorname*{Rep}(\mathcal{O}_{N},L^{2}(\mathbb{T}))
\]
where $\mathbb{T}=\mathbb{R}/2\pi\mathbb{Z}$ is the $1$-torus. We review joint
research papers \cite{BEJ00}, \cite{BrJo97b}, \cite{BrJo98a}, \cite{BrJo99b},
\cite{BrJo00} and the solo paper \cite{Jor00}. See also \cite{Dau92}.

If $A$ is given, we show that the operators
\[
S_{i}^{(A)}=\pi^{(A)}(s_{i})
\]
are weighted shift operators on $L^{2}(\mathbb{T})$. Our papers \cite{BrJo97b}
and \cite{BJKW00} indicate generalizations to $\mathbb{T}^{d}$, $d>1$, but we
restrict to $d=1$ here.

\section{A Hilbert Module\label{section2}}

Let $(X,\mu)$ be a measure space, and we assume that $\mu$ is a probability
measure on $X$. Let $\sigma:X\rightarrow X$ be a measurable $N$-to-$1$ self
map, and let $N$ be given, and fixed.

It will further be assumed that there is a probability measure $\mu$ on $X$
such that
\begin{equation}
\frac{1}{N}\sum_{i=0}^{N-1}\mu\circ\sigma_{i}^{-1}=\mu, \label{eq2.2}%
\end{equation}
where $\sigma_{i}:X\rightarrow X$ is some chosen sections for $\sigma$, i.e.,
satisfying
\begin{equation}
\sigma\circ\sigma_{i}=\operatorname*{id}\nolimits_{X},\qquad i=0,\dots
,N-1\text{. } \label{eq2.3}%
\end{equation}
The analysis below refers to such a measure $\mu$. If $\sigma$ is expansive,
$\mu$ is known to exist and be unique \cite{BJO99}. If $X=\mathbb{T}$, and
$\sigma(z)=z^{N}$, then $\mu$ is the usual Haar measure on $\mathbb{T}$.

Let
\begin{equation}
\mathcal{A}_{\sigma}:=\{f\circ\sigma\,;\,f\in L^{\infty}(X)\}. \label{eq2.1}%
\end{equation}
We will also assume that $L^{\infty}(X)$ and its subalgebras act by
multiplication on the Hilbert space $\mathcal{H}_{\mu}=L^{2}(X,\mu)$.

\begin{lemma}
\label{lemma2.1}

\begin{enumerate}
\item[$\mathrm{(a)}$] A system of measurable functions
\[
m_{i}:X\longrightarrow\mathbb{C},\qquad i=0,\dots,N-1,
\]
forms an orthonormal basis for $\mathcal{H}_{\mu}$ as an $\mathcal{A}_{\sigma
}$-Hilbert module if and only if there are sections $\sigma_{i}:X\rightarrow
X$, i.e., $\sigma\circ\sigma_{i}=\operatorname*{id}_{X}$ such that the
$N\times N$ matrix
\begin{equation}
M_{m}:=(m_{i}\circ\sigma_{j})_{i,j=0}^{N-1}\qquad\text{is unitary,}%
\label{eq2.4}%
\end{equation}
i.e., defines
\[
M_{m}:X\longrightarrow\mathrm{U}_{N}(\mathbb{C}).
\]

\item[$\mathrm{(b)}$] If $m_{0}\in L^{\infty}(X)$ is given such that
\begin{equation}
\sum_{y:\sigma(y)=x}\left|  m_{0}(y)\right|  ^{2}\equiv1\text{ a.a. }x\in X,
\label{eq2.5}%
\end{equation}
then there is a measurable selection $m_{1},\dots,m_{N-1}$ such that the
combined system $m_{0},\dots,m_{N-1}$ satisfies \textup{(\ref{eq2.4}).}

\item[$\mathrm{(c)}$] A function system $\{m_{i}\}_{i=0}^{N-1}$ satisfies the
conditions in \textup{(a)} if and only if the system of operators
\begin{equation}
S_{i}^{(m)}f:=\sqrt{N}m_{i}\,f\circ\sigma,\qquad f\in\mathcal{H}_{\mu
},\label{eq2.6}%
\end{equation}
defines a representation of $\mathcal{O}_{N}$, i.e., $S^{(m)}\in
\operatorname*{Rep}(\mathcal{O}_{N},\mathcal{H}_{\mu})$. \textup{(}Note the
right-hand side in \textup{(\ref{eq2.6})} is the function $\sqrt{N}%
m_{i}(x)f(\sigma(x))$.\textup{)}
\end{enumerate}
\end{lemma}

\begin{definition}
\label{definition2.2}The system $\{m_{i}\}_{i=0}^{N-1}$ is called a quadrature
mirror filter system \textup{(}QMF\/\textup{)} by analogy to the example
$N=2$, in which case $m_{0}\sim(a_{k})$ serves as the low pass filter, and
$m_{1}\sim(b_{k})$ as the high pass filter. It is the orthogonality relations
\textup{(\ref{eq2.4})} which motivate the name QMF, and the use of filters in
wavelet theory is further fleshed out in \cite{Dau92}.
\end{definition}

\begin{proof}
[Proof of Lemma 2.1]Most of the details are contained in the paper
\cite{BrJo97b}, and others will be in a later more detailed paper. But we
sketch here the argument for orthogonality in the Hilbert module. Let the
functions $\{m_{i}\}_{i=0}^{N-1}$ be as stated in the lemma. Orthogonality
refers to
\[
\left\langle \mathcal{A}_{\sigma}m_{i}\mid\mathcal{A}_{\sigma}m_{j}%
\right\rangle _{\mathcal{H}_{\mu}}=0
\]
whenever $i\neq j$. Hence we must calculate, for $f\in L^{\infty}(X)$:%
\[
\int_{X}\overline{m_{i}(x)}f(\sigma(x))m_{j}(x)\,d\mu(x)=\frac{1}{N}\int
_{X}f(x)\sum_{y:\sigma(y)=x}\overline{m_{i}(y)}m_{j}(y)\,d\mu(x),
\]
using (\ref{eq2.2}), and the sum under the integral vanishes by (\ref{eq2.4})
if $i\neq j$.

It is important that the selection result mentioned in (b) is generally not
possible in the category of continuous functions. (See \cite{Kad84}.)

Conversely, if $\left\{  m_{i}\right\}  _{i=0}^{N-1}\subset L^{\infty}\left(
X\right)  $ is given such that (\ref{eq2.6}) defines a representation of
$\mathcal{O}_{N}$ on $\mathcal{H}_{\mu}$, then the closed subspaces
$S_{i}^{\left(  m\right)  }\mathcal{H}_{\mu}$ ($\subset\mathcal{H}_{\mu}$) are
the submodules in a corresponding orthogonal $\mathcal{A}_{\sigma}$-module
decomposition, $\mathcal{H}_{\mu}=\sum_{i=0}^{N-1}\left[  \mathcal{A}_{\sigma
}m_{i}\right]  $, i.e., the $L^{2}\left(  \mu\right)  $-closure of
$m_{i}\mathcal{A}_{\sigma}$ is $S_{i}^{\left(  m\right)  }\mathcal{H}_{\mu}$
for each $i$. Specifically, from (\ref{eq1.1}), we get the identity
$\mathcal{H}_{\mu}\rightarrow f=\sum_{i}S_{i}^{\left(  m\right)  }%
S_{i}^{\left(  m\right)  \,\ast}f=\sum_{i}m_{i}k_{i}\circ\sigma$, where
$k_{i}=\sqrt{N}S_{i}^{\left(  m\right)  \,\ast}f$.
\end{proof}

\section{\label{Com}Comparing Representations}

\begin{lemma}
\label{LemCom.1}Let $G_{N}\left(  X\right)  $ be the group of measurable maps%
\[
A\colon X\longrightarrow\mathrm{U}_{N}\left(  \mathbb{C}\right)  .
\]
Then $G_{N}\left(  X\right)  $ acts transitively on the systems $\left\{
m_{i}\right\}  _{i=0}^{N-1}$ of functions from Lemma \textup{\ref{lemma2.1}.}
\end{lemma}

\begin{proof}
Let $\left\{  m_{i}\right\}  $ be given as in Lemma \textup{\ref{lemma2.1}},
and let $A\in G_{N}\left(  X\right)  $. Set%
\begin{equation}
n_{i}\left(  x\right)  :=\sum_{j=0}^{N-1}A_{i,j}\left(  \sigma\left(
x\right)  \right)  m_{j}\left(  x\right)  . \label{eqCom.1}%
\end{equation}
Then $\left\{  n_{i}\right\}  $ satisfies the same orthogonality relations.
For
\begin{align*}
\sum_{k}\overline{n_{i}\left(  \sigma_{k}\left(  x\right)  \right)  }%
n_{j}\left(  \sigma_{k}\left(  x\right)  \right)   &  =\sum_{k}\sum
_{l,l^{\prime}}\bar{A}_{i,l}(x)\bar{m}_{l}\left(  \sigma_{k}\left(  x\right)
\right)  A_{j,l^{\prime}}\left(  x\right)  m_{l^{\prime}}\left(  \sigma
_{k}\left(  x\right)  \right) \\
&  =\sum_{l,l^{\prime}}\delta_{l,l^{\prime}}\bar{A}_{i,l}\left(  x\right)
A_{j,l^{\prime}}\left(  x\right)  =\sum_{l}\bar{A}_{i,l}\left(  x\right)
A_{j,l}\left(  x\right)  =\delta_{i,j},
\end{align*}
proving the assertion. We used the fact that $\sigma(\sigma_{k}(x))=x$, see
(\ref{eq2.3}).

Conversely, if $S_{i}^{\left(  m\right)  }$ and $S_{i}^{\left(  n\right)  }$
are given representations as described in (\ref{eq2.6}) of Lemma
\ref{lemma2.1}, then%
\begin{equation}
\left(  S_{i}^{\left(  n\right)  \,\ast}S_{j}^{\left(  m\right)  }\right)  \in
M_{N}\left(  L^{\infty}\left(  X\right)  \right)  \qquad(=M_{N}\left(
\mathbb{C}\right)  \otimes L^{\infty}\left(  X\right)  ). \label{eqCom.2}%
\end{equation}
For it follows from the Cuntz relations (\ref{eq1.1}) that the matrix in
(\ref{eqCom.2}) is unitary, and a computation shows that its matrix entries
are multiplication operators. Indeed,%
\[
\left(  S_{i}^{\left(  m\right)  }\right)  ^{\ast}f\left(  x\right)  =\frac
{1}{\sqrt{N}}\sum_{y\colon\sigma\left(  y\right)  =x}\overline{m_{i}\left(
y\right)  }f\left(  y\right)  \text{,\qquad for }x\in X\text{ and for }%
f\in\mathcal{H}_{\mu}.
\]
Hence
\[
\left(  S_{i}^{\left(  n\right)  }\right)  ^{\ast}S_{j}^{\left(  m\right)
}f\left(  x\right)  =\sum_{y\colon\sigma\left(  y\right)  =x}\overline
{n_{i}\left(  y\right)  }m_{j}\left(  y\right)  f\left(  x\right)  ,
\]
and so%
\[
A_{i,j}\left(  x\right)  :=\sum_{y\colon\sigma\left(  y\right)  =x}%
\overline{n_{i}\left(  y\right)  }m_{j}\left(  y\right)
\]
defines an element of $G_{N}\left(  X\right)  $, and it satisfies
(\ref{eqCom.1}) by its very construction. This proves transitivity.
\end{proof}

When the lemma is applied to the example $X=\mathbb{T}$, $\sigma\left(
x\right)  =z^{N}$, $z\in\mathbb{T}$, and $\sigma_{k}\left(  x\right)
=e^{i\frac{2\pi k}{N}}\sqrt[N]{z}$, $k=0,\dots,N-1$, we conclude that the loop
group (see Section \ref{section1}) of measurable $A\colon\mathbb{T}%
\rightarrow\mathrm{U}_{N}\left(  \mathbb{C}\right)  $ acts transitively on the
wavelet representations. Let $n_{k}\left(  z\right)  =\frac{1}{\sqrt{N}}z^{k}%
$, $k=0,\dots,N-1$, and let $\left\{  m_{k}\right\}  _{k=0}^{N-1}$ be an
arbitrary $m$-system as in Lemma \ref{lemma2.1}.

Then we have

\begin{corollary}
\label{CorCom.2}There is a $1$-to-$1$ bijective correspondence between the
loops $A$ and the $m$-systems of Lemma \textup{\ref{lemma2.1}} given as
follows:%
\[
m_{i}\left(  z\right)  =\frac{1}{\sqrt{N}}\sum_{j=0}^{N-1}A_{i,j}\left(
z^{N}\right)  z^{j}%
\]
and%
\[
A_{i,j}\left(  z\right)  =\frac{1}{\sqrt{N}}\sum_{w\in\mathbb{T}\colon
w^{N}=z}m_{i}\left(  w\right)  w^{-j}.
\]
\end{corollary}

\section{Wavelets\label{section4}}

Let $\varphi\in L^{2}(\mathbb{R})$ be the compactly supported scaling
function, i.e., a solution to the scaling identity
\begin{equation}
\varphi(x)=\sum_{k=0}^{Ng-1}a_{k}\varphi(Nx-k). \label{eq2.8}%
\end{equation}
Then the wavelet generators $\psi_{1},\dots,\psi_{N-1}\in L^{2}(\mathbb{R})$
are constructed from $\varphi$ by use of Lemma \ref{lemma2.1}(b) above, and
standard wavelet tools from \cite{Dau92}. The generators make the system
\[
\{N^{j/2}\psi_{i}(N^{j}x-k)\,;\,i=0,1,\dots,N-1,\;j,k\in\mathbb{Z}\}
\]
into an orthonormal basis for $L^{2}(\mathbb{R})$; except for a smaller
variety of cases when the system is only a tight frame. The coefficients
$\{a_{k}\}$ represent a wavelet filter, and
\begin{equation}
m_{0}(z)=\sum_{k}a_{k}z^{k}. \label{new}%
\end{equation}
Then define the operator
\begin{equation}
W\colon\ell^{2}(\mathbb{Z})\longrightarrow L^{2}(\mathbb{R}) \label{eq2.9}%
\end{equation}
by
\begin{equation}
W\xi=\sum_{k}\xi_{k}\varphi(x-k). \label{eq2.10}%
\end{equation}
The conditions on the wavelet filter $\{a_{k}\}$\textbf{ }may now be restated
in terms of $m_{0}(z)$ in (\ref{new}) as follows:
\begin{equation}
\sum_{k=0}^{N-1}\left|  m_{0}(ze^{i\frac{k2\pi}{N}})\right|  ^{2}=1,\text{ }
\label{eq2.11}%
\end{equation}
and%
\begin{equation}
m_{0}(1)=1\text{,\qquad the low pass property.} \label{eq2.12}%
\end{equation}
Then $W$ in (\ref{eq2.10}) maps $\ell^{2}(\mathbb{Z})$ onto the resolution
subspace $\mathcal{V}_{0}$ ($\subset L^{2}(\mathbb{R})$), and we note that
\begin{equation}
U_{N}W=WS_{0}, \label{eq2.13}%
\end{equation}
where
\begin{equation}
U_{N}f(x)=N^{-\frac{1}{2}}f\left(  \frac{x}{N}\right)  ,\qquad f\in
L^{2}(\mathbb{R}),\;x\in\mathbb{R}. \label{eq.2.14}%
\end{equation}
We showed in \cite{BrJo00} that there are functions $m_{1},\dots,m_{N-1}$ such
that the $N$-by-$N$ complex matrix
\begin{equation}
\left(  m_{j}(e^{i\frac{k2\pi}{N}}z)\right)  _{j,k=0}^{N-1} \label{eq2.15}%
\end{equation}
is unitary for all $z\in\mathbb{T}$. (See Lemma \ref{lemma2.1}(b).) We define
\begin{equation}
S_{j}f(z)=\sqrt{N}m_{j}(z)f(z^{N}),\qquad f\in L^{2}(\mathbb{T}),\;z\in
\mathbb{T}\text{. } \label{eq2.16}%
\end{equation}

The main result will be stated in the present section, but without proof.
Instead the reader is referred to \cite{Jor00} for the full proof, and for a
detailed discussion of its implications. We noted above that the
representation (\ref{eq2.16}) given from a QMF system $m_{j}=m_{j}^{\left(
A\right)  }$, via
\begin{equation}
m_{j}^{(A)}(z)=\frac{1}{\sqrt{N}}\sum_{k=0}^{N-1}A_{j,k}(z^{N})z^{k},
\label{eq61}%
\end{equation}%
\begin{equation}
A_{j,k}(z)=\frac{1}{\sqrt{N}}\sum_{w^{N}=z}w^{-k}m_{j}(w) \label{eq62}%
\end{equation}
is irreducible if and only if the subbands are optimal, in that they do not
admit further reduction into a refined system of closed subspaces of
$L^{2}\left(  \mathbb{R}\right)  $.

\begin{theorem}
\label{ThmThe}The representation
\[
S_{j}^{\left(  A\right)  }f\left(  z\right)  =\sqrt{N}m_{j}^{\left(  A\right)
}\left(  z\right)  f\left(  z^{N}\right)  ,\qquad f\in L^{2}\left(
\mathbb{T}\right)  ,\;z\in\mathbb{T},
\]
is an irreducible representation of $\mathcal{O}_{N}$ on $L^{2}\left(
\mathbb{T}\right)  $ if and only if $A\colon\mathbb{T}\rightarrow
\mathrm{U}_{N}\left(  \mathbb{C}\right)  $ does not admit a matrix corner of
the form
\begin{equation}
V%
\begin{pmatrix}
z^{n_{0}} &  &  & \\
& z^{n_{1}} &  & \llap{\smash{\huge$0$}}\\
&  & \ddots & \\
\rlap{\smash{\huge$0$}} &  &  & z^{n_{M-1}}%
\end{pmatrix}
, \label{eqThmMin.7AzVm}%
\end{equation}
for some $V\in\mathrm{U}_{M}\left(  \mathbb{C}\right)  $, and where
$n_{0},n_{1},\dots,n_{M-1}\in\left\{  0,1,2,\dots\right\}  $. Moreover, two
representations $\pi^{\left(  A\right)  }$ and $\pi^{\left(  B\right)  }$
defined from different loops $A$, $B$ are unitarily inequivalent unless
$A\equiv B$ modulo a matrix corner of type \textup{(\ref{eqThmMin.7AzVm}).}
\end{theorem}

\begin{acknowledgements}
We are grateful to Brian Treadway and Cymie Wehr for excellent typesetting,
and to the ICMP participants for enlightening discussions. The author also had
helpful discussions with W. Arveson, O. Bratteli, and P. Muhly.
\end{acknowledgements}

%\bibliographystyle{BFTALPHA}
%\bibliography{jorgen}

\end{document}